\documentclass[11pt]{article} \usepackage[dvips]{epsfig}
\usepackage{latexsym}
\usepackage{amsfonts,amsmath,amssymb}

\newtheorem{lemma}{Lemma}
\newtheorem{definition}{Definition}

\newtheorem{theorem}{Theorem}

\newcommand{\qed}{\mbox{$\Diamond$}\vspace{\baselineskip}}
\newenvironment{proof}{\noindent{\bf Proof:}}{\qed}
\newcommand{\Var}{\hbox {Var}}
\begin{document}
\author{Mikl\'os B\'ona\\
        Department of Mathematics\\
University of Florida\\
Gainesville FL 32611-8105
  }

\title{Generalized Descents and Normality}
\date{today}

\maketitle

\begin{abstract} We use Janson's dependency criterion to prove that 
the distribution of $d$-descents of permutations of length $n$
converge to a normal distribution as $n$ goes to infinity. We show
that this remains true even if $d$ is allowed to grow with $n$, up to 
a certain degree. 
\end{abstract}

\section{Introduction}
Let $p=p_1p_2\cdots p_n$ be a permutation. We say that the pair $(i,j)$ is
a $d$-descent in $p$ if $i<j\leq i+d$, and $p_i>p_j$. In particular,
1-descents correspond to descents in the traditional sense, and
$(n-1)$-descents correspond to inversions. This concept was introduced
in \cite{demari} by De Mari and Shayman, whose motivation came from
algebraic geometry. They have proved that if $n$ and $d$ are fixed, and
$c_k$ denotes the number of permutations of length $n$ with exactly $k$
$d$-descents, then the sequence $c_0,c_1, \cdots $ is unimodal, that is,
it increases steadily, then it decreases steadily. It is not known in general
if
the sequence  $c_0,c_1, \cdots $ is log-concave or not, that is, whether
$c_{k-1}c_{k+1} \leq c_k^2$ holds for all $k$. We point out that in general,
the polynomial $\sum_kc_kx^k$ does not have real roots only. Indeed, in the
special case of $d=n-1$, we get the well-known \cite{combperm} identity
\[\sum_kc_kx^k=(1+x)\cdot (1+x+x^2)\cdot \cdots\cdot  (1+x+\cdots +x^{n-1}),\]
which has all $n$th roots of unity as roots. Indeed, in this case, a 
$d$-descent is just an inversion, as we said above.

In this paper, we prove a related property of generalized descents by
showing that their distribution converges to a normal distribution as
the length $n$ of our permutations goes to infinity. Our main tool is
{\em Janson's dependency criterion}, which is a tool to prove normality for
sums of bounded random variables with a sparse  dependency graph. 

\section{The Proof of Asymptotic Normality}
\subsection{Background and Definitions} 

We need to introduce some notation for transforms of the random variable
$Z$. Let $\bar{Z}=Z-E(Z)$, let $\tilde{Z}=\bar{Z}/\sqrt{\Var( Z)}$, and let
$Z_n\rightarrow N(0,1)$ mean that $Z_n$ converges in distribution to the 
standard normal variable. 

For the rest of this paper, let $d\geq 1$ be a fixed positive integer. 
Let $X_n=X_n^{(d)}$ denote the random variable counting the $d$-descents of
a randomly selected permutation of length $n$. We want to prove that
$X_n$ converges to a normal distribution as $n$ goes to infinity, in
other words, that $\tilde{X}_n \rightarrow N(0,1)$ as $n\rightarrow \infty$.
 Our main 
tool in doing so is a  theorem called Janson's dependency 
criterion. 
In order to state that theorem, we need the following definition.

\begin{definition}
Let $\{Y_{n,k}|k=1,2\cdots \}$ be an array of random variables.
 We say that a graph $G$ is 
a {\em dependency graph} for   $\{Y_{n,k}|k=1,2\cdots \}$ if the following
two conditions are satisfied:
\begin{enumerate}
\item There exists a bijection between the random variables $Y_{n,k}$ and
the vertices of $G$, and
\item If $V_1$ and $V_2$ are two disjoint sets of vertices of $G$ so that
no edge of $G$ has one endpoint in $V_1$ and another one in $V_2$, then
the corresponding sets of random variables are independent.
\end{enumerate}
\end{definition}

Note that  the dependency graph of a
 family of variables is not unique. Indeed if $G$ is a dependency graph
for a family and $G$ is not a complete graph,
 then we can get other dependency graphs for the family
by simply adding new edges to $G$. 

Now we are in position to state Janson's dependency criterion.

\begin{theorem} \cite{janson} \label{janson}
Let $Y_{n,k}$ be an array of random variables such that for all $n$, and
for all $k=1,2,\cdots ,N_n$, the inequality $|Y_{n,k}|\leq A_n$ holds for
some real number $A_n$, and that the maximum degree of a dependency
graph of $\{Y_{n,k} | k=1,2,\cdots ,N_n \}$ is $\Delta_n$. 

Set $Y_n=\sum_{k=1}^{N_n} Y_{n,k}$ and $\sigma_n^2= \Var ( Y_n)$. If there
is a natural number $m$ so that
\begin{equation} \label{jansencond}
N_n\Delta_n^{m-1} \left (\frac{A_n}{\sigma_n} \right )^m \rightarrow 0,
\end{equation}
then \[ \tilde{Y}_n \rightarrow N(0,1) .\]
\end{theorem}

\subsection{Applying Janson's Criterion}
We will apply Janson's theorem with the $Y_{n,k}$ being the indicator 
random variables $X_{n,k}$ of the event that a given ordered pair of indices 
(indexed by $k$
in some way) 
form a $d$-descent in the randomly selected permutation $p=p_1p_2\cdots p_n$.
 So $N_n$ is the number of pairs $(i,j)$ of
indices so that $1\leq i<j\leq i+d  \leq n$. Then by definition,
\[Y_n=\sum_{k=1}^{N_n}Y_{n,k}=\sum_{k=1}^{N_n}X_{n,k}=X_n.\]

There remains the task of verifying that the variables $Y_{n,k}$ satisfy
all conditions of Jansen's theorem. 

First, it is clear that
$N_n\leq nd$, and we will compute the exact value of $N_n$ later.
 By the definition of indicator random variables, we have
 $|Y_{n,k}|\leq 1$, so we can set $A_n=1$ for all $n$. 

Next we consider the numbers $\Delta_n$ in the following 
dependency graph of the family of the $Y_{n,k}$. 
 Clearly, the indicator random
variables that belong to two pairs $(i,j)$ and $(r,s)$ of indices are
independent if and only if the sets $\{i,j\}$ and $\{r,s\}$ are 
disjoint.  So fixing $(i,j)$, we need one of $i=r$, $i=s$, $j=r$ or $j=s$
to be true for the two distinct variables to be dependent. So let the
vertices of $G$ be the $N_n$ pairs of indices $(i,j)$ so that
$i<j\leq i+d$, and connect $(i,j)$ to $(r,s)$ if one of 
$i=r$, $i=s$, $j=r$ or $j=s$ holds. The graph defined in this way is clearly
a dependency graph for the family of the $Y_{n,k}$. For a fixed
pair $(i,j)$, each of these four
equalities occurs at most $d$ times. (For instance, if $i=s$, then 
$r$ has to be one of $i-1,i-2,\cdots ,i-d$.)
Therefore, $\Delta_n\leq 4d$.

If we take a new look at  (\ref{jansencond}), we see that the Janson
criterion will be satisfied if we can show that $\sigma_n$ is large. 
This is the content of the next lemma.

\begin{lemma}
If $n\geq 2d$, then
\begin{equation} \label{precise}
\Var (X_n) = \frac{6dn+10d^3-3d^2-d}{72}.  \end{equation}
In particular, $ \Var (X_n)$ is a {\em linear} function of $n$.
\end{lemma}

Note that in particular, for $d=1$, we get the well-known fact \cite{combperm}
that the variance of Eulerian numbers in permutations of length $n$ is
$(n+1)/12$.

\begin{proof}
By linearity of expectation, we have
\begin{eqnarray} \label{variance}
\Var (X_n) & = & E(X_n^2) - (E(X_n))^2 \\
 & = & E \left (\left( \sum_{k=1}^{N_n} X_{n,k} \right )^2 \right )
-  \left (E \left (\sum_{k=1}^{N_n} X_{n,k} \right ) \right )^2 \\
 & = & E \left (\left( \sum_{k=1}^{N_n} X_{n,k} \right )^2 \right )
- \left( \sum_{k=1}^{N_n} E(X_{n,k}) \right )^2  \\
 & = &  \sum_{k_1, k_2}
E(X_{n,k_1}X_{n,k_2})  - \sum_{k_1, k_2}
E(X_{n,k_1})E(X_{n,k_2})
\end{eqnarray}

Clearly, $E(X_{n,k})=1/2$, so the $N_n^2$ summands that appear in the last
line of the above chain of equations with a {\em negative sign}
 are each equal to $1/4$. As far
as the $N_n^2$ summands that appear with a positive sign, {\em most} of
them are equal to $1/4$. More precisely, if $X_{n,k_1}$ and $X_{n,k_2}$ 
are independent, then 
\[E(X_{n,k_1}X_{n,k_2})=E(X_{n,k_1})E(X_{n,k_2})=\frac{1}{4}.\]
If $k_1=k_2$, then $E(X_{n,k_1}X_{n,k_2})=E(X_{k_1}^2=E(X_{k_1})=1/2$.
Otherwise, if $X_{n,k_1}$ and $X_{n,k_2}$ 
are dependent, then either $E(X_{n,k_1}X_{n,k_2})=1/3$, or
 $E(X_{n,k_1}X_{n,k_2})=1/6$. Indeed, if $X_{k_1}$ is the indicator variable
of the pair $(i,j)$ being a $d$-descent and $X_{k_2}$
is the indicator variable
of the pair $(r,s)$ being a $d$-descent, then as we said above, 
$X_{n,k_1}$ and $X_{n,k_2}$ 
are dependent if and only if one of $i=r$, $i=s$, $j=r$ or $j=s$ holds.
 If $i=r$ or
$j=s$ holds, then $E(X_{n,k_1}X_{n,k_2})=1/3$, and if $i=s$ or $j=r$ holds,
then  $E(X_{n,k_1}X_{n,k_2})=1/6$. Indeed, for instance, with $i=r$, 
we have $X_{n,k_1}=X_{n,k_2}=1$ if and only if $p_i$ is the largest of
the entries $p_i$, $p_j$, and $p_s$. Similarly, with $i=s$, we have
 $X_{n,k_1}=X_{n,k_2}=1$ if and only if $p_r>p_i>p_j$.

We will now count  how many summands $E(X_{n,k_1}X_{n,k_2})$ are equal to
$1/2$, to $1/3$, and to $1/6$. 

\begin{enumerate}
\item First, $E(X_{n,k_1}X_{n,k_2})=1/2$ if and 
only if $k_1=k_2$. This happens $N_n$ times, once for each pair $(i,j)$
so that $i<j\leq i+d$. For a given $i$, there are $d$ such pairs 
if $i\leq n-d$, and $d-t$ such pairs if $i=n-d+t$, so
\[N_n=(n-d)d + (d-1) + (d-2)+\cdots + 1=(n-d)d +{d\choose 2}.\]
\item 
Second,  $E(X_{n,k_1}X_{n,k_2})=1/3$ if $i=r$, or $j=s$. By symmetry, we
can consider the first case, then multiply by two. If $i\leq n-d$, then
we have $d(d-1)$ choices for $j$ and $s$, and if $i=n-d+t$, then we
have $(d-t)(d-t-1)$ choices. So the number of pairs $(k_1,k_2)$ so that
 $E(X_{n,k_1}X_{n,k_2})=1/3$ is
\noindent \[2(n-d)d(d-1)+2(d-1)(d-2)+2(d-2)(d-3)+\cdots +2\cdot 2\cdot 1=\]
\[2(n-d)d(d-1)+4{d\choose 3}.\]
\item 
Finally,  $E(X_{n,k_1}X_{n,k_2})=1/6$ if $i=s$, or $j=r$. By symmetry, we
can again consider the first case, then multiply by two. If $d\leq i\leq n-d$,
then there are $d^2$ choices for $(j,r)$. If $i\leq d$, then there are
$d$ choices for $j$, and $i-1$ choices for $r$. If $n-d<i$, then there are
$n-i$ choices for $j$, and $d$ choices for $r$, assuming that $n\geq 2d$.  
So the number of pairs $(k_1,k_2)$ so that
 $E(X_{n,k_1}X_{n,k_2})=1/6$ is
\[2(n-2d)d^2+2(d-1)d+2(d-2)d+\cdots +2d=2(n-2d)d^2 + d^2(d-1).\]
For all remaining pairs $(k_1,k_2)$, the variables $X_{n,k_1}$ and $X_{n,k_2}$
are independent, and so $E(X_{n,k_1}X_{n,k_2})=1/4$.
\end{enumerate}

Comparing our results from cases 1-3 above with (\ref{variance}),
 and recalling that in all other cases,
$E(X_{n,k_1}X_{n,k_2})=1/4$, we obtain the formula that was to be proved.  
\end{proof}

The proof of our main theorem is now immediate.

\begin{theorem} \label{basic}
Let $d$ be a fixed positive integer.
Let $X_n$ be the random variable counting  $d$-descents of a
randomly selected $n$-permutation.
Then $\tilde{X}_n \rightarrow N(0,1)$.
\end{theorem}

\begin{proof}
Use Theorem \ref{janson} with $Y_n=X_n$, $\Delta_n=4d$, 
$N_n=(n-d)d+{d\choose 2}$, and $\sigma_n=\sqrt{\frac{6dn+10d^3-3d^2-d}{72}}$.
All we need to show is that there exists a positive integer $m$ so that
\[\left((n-d)d+{d\choose 2}\right)\cdot (4d)^{m-1} \cdot 
\left ( \frac{72}{6dn+10d^3-3d^2-d} \right )^{m/2} \rightarrow 0,\]
for which it suffices to find a positive integer $m$ so that
\begin{equation} \label{easier} 
(dn)\cdot (4d)^{m-1}\cdot \left (\frac{12}{dn} \right )^{m/2}
 \rightarrow 0 .\end{equation}
Clearly, any $m\geq 3$ suffices, since for any such $m$, the left-hand side
is of the form $C/n^{\alpha}$, for positive constants $C$ and $\alpha$. 
\end{proof}

\section{Further Directions}
We see from (\ref{easier}) that the statement of Theorem \ref{basic} can
be strengthened, from a constant $d$ to a $d$ that is a function of $n$.
 Indeed, (\ref{easier}) is equivalent to saying
that \[cn\left(\frac{d}{n}\right ) ^{m/2} \rightarrow 0.\]
This convergence holds as long as $d\leq n^{1-\epsilon}$ for some fixed
positive $\epsilon$, we can choose $m$ so that $(m/2)\cdot \epsilon >1$, and
then condition (\ref{easier}) will be satisfied. So we have proved the
following theorem.

\begin{theorem} Let $n\rightarrow \infty$, and let us assume
that there exists a positive constant $\epsilon$ so that for $n$ sufficiently
large,  $d=d(n)\leq n^{1-\epsilon}$. Let $X_n$ be defined as before. 
Then  \[\tilde{X}_n \rightarrow N(0,1).\]
\end{theorem}

This leaves the cases of larger $d$ open. We point out that in the special
 case of $d=n-1$, that is, inversions, asymptotic normality is known 
\cite{diaconis}, \cite{fulman}.

Another possible direction for generalizations is the following. 
Let ${\bf d}
= (d_1,d_2 \cdots, d_{n-1})$, where the $d_i$ are positive
integers. If $p = p_1 ... p_n$ is
in an $n$-permutation, let $f_d(p)$ be the number of pairs $(i,j)$
 such that $0 < j-i \leq
d_i$ and $p_i > p_j$. For instance, if ${\bf d}=(1,1,...,1)$ then $f_d(p)$
 is the
number of descents of $p$. If ${\bf d}=(n-1,n-2,...,1)$ then $f_d(p)$ is the
number of inversions of $p$. It is known \cite{demari}, by an 
argument from algebraic geometry, that if
   \[c_k = |\{ p\in S_n : f_d(p)=k \}|,\]
then the sequence $c_0,c_1,\cdots $ is unimodal. Log-concavity and normality
are not known. Note that in this paper, we have treated the special case of
${\bf d}=(d,d,\cdots ,d)$. 

\vskip 2 cm 
\centerline { {\bf Acknowledgment}}

I am thankful to Richard Stanley who introduced me to the topic of generalized
descents.

\end{document}